\documentclass[11pt]{article}%
\usepackage{amsfonts}
\usepackage[applemac]{inputenc}
\usepackage{amssymb}
\usepackage{makeidx}
\usepackage{graphicx}
\usepackage{amsmath}
\usepackage[a4paper]{geometry}
\usepackage[displaymath]{lineno}
\usepackage[singlespacing]{setspace}
\usepackage[normalem]{ulem}
\usepackage{lineno}
\usepackage{color}%
\providecommand{\U}[1]{\protect\rule{.1in}{.1in}}
{\setlength\paperheight {297mm}
\setlength\paperwidth  {210mm}}

\providecommand{\U}[1]{\protect\rule{.1in}{.1in}}

\newtheorem{theo}{Theorem}
\newtheorem{lem}[theo]{Lemma}

\newtheorem{cor}[theo]{Corollary}
\newtheorem{rem}{Remark}
\newtheorem{exam}{Example}

\newenvironment{dem}[1][Proof]{\noindent \textbf{#1.} }{\ \rule{0.5em}{0.5em}}
\begin{document}

\title{Alternative representations of the normal cone to the domain of supremum
functions and subdifferential calculus\thanks{Research supported by ANID
(Fondecyt 1190012 and 1190110), Proyecto CMM ANID PIA AFB170001, MICIU of
Spain and Universidad de Alicante (Contract Beatriz Galindo BEA-GAL 18/00205),
and Research Project PGC2018-097960-B-C21 from MICINN, Spain. The research of
the third author is also supported by the Australian ARC - Discovery Projects
DP 180100602.} }
\author{R. Correa\thanks{e-mail: rcorrea@dim.uchile.cl}, A.\ Hantoute$^{§}%
$\thanks{e-mail: hantoute@ua.es (corresponding author)} \ and M.A. López$^{¶}%
$\thanks{e-mail: marco.antonio@ua.es}\\${^{\dag}}${\scriptsize Universidad de O'Higgins, Chile, and DIM-CMM of
Universidad de Chile}\\$^{\ddag}${\scriptsize Center for Mathematical Modeling (AFB170001),
Universidad de Chile }\\$^{§}${\scriptsize Universidad de Alicante, Spain}\\$^{¶}${\scriptsize CIAO, Federation University, Ballarat, Australia}}
\date{}
\maketitle

\begin{abstract}
The first part of the paper provides new characterizations of the normal cone
to the effective domain of the supremum of an arbitrary family of convex
functions. These results are applied in the second part to give new
formulas\ for\ the subdifferential of the supremum function, which
use\ both\ the active and nonactive functions at the reference point. Only the
data functions are involved\ in these characterizations, the active ones from
one side, together with\ the\emph{ }nonactive functions multiplied by some
appropriate parameters. In contrast with previous works in the literature, the
main feature of our subdifferential characterization is that the normal cone
to the effective domain of the supremum (or to finite-dimensional sections of
this domain) does not appear. A new type\ of optimality conditions for convex
optimization is established at the end of the paper.



\end{abstract}

\section{Introduction}

Given the pointwise supremum $f:=\sup_{t\in T}f_{t}$ of a family\ of convex
functions $f_{t}:X\rightarrow\mathbb{R}\cup\{+\infty\}$, $t\in T$, $T$ being a
non-empty and arbitrary, defined on a separated locally convex space $X,$ many
researchers have addressed the paradigmatic problem of characterizing the
subdifferential of the supremum, $\partial f(x)$, at any point $x$ of the
effective domain of $f.$ These characterizations are usually given in terms of
the (approximate-) subdifferentials of the data functions, $\partial
_{\varepsilon}f_{t}(x),\ t\in T,$ $\varepsilon\geq0,$ and, in the most general
cases, in terms also of the normal cone to the effective domain of $f$ or to
finite-dimensional sections of it. The interest of this problem comes from the
fact that\textbf{ }many convex functions, such as the Fenchel conjugate, the
sum, the composition with affine mappings, etc., can be expressed as the
supremum of affine or convex functions. Therefore, getting formulas for the
subdifferential of the supremum is expected to play a crucial role in convex
optimization and variational analysis. Some remarkable contributions to the
topic are: Brøndsted \cite{Br72}, Ioffe \cite{Io12}, Ioffe \& Levin
\cite{IL72}, Ioffe \& Tikhomirov \cite{IT79}, Levin \cite{Le69}, Pschenichnyi
\cite{Ps65}, Rockafellar \cite{Ro79}, Valadier \cite{Va69}, etc. In
\cite{Ti90} the historical origins of the issue are traced out. More recently,
in a series of papers (\cite{CHL16}, \cite{CHL19b}, \cite{HL08}, \cite{HLZ08},
etc.) new characterizations of the subdifferential supremum in different
settings are provided, and some related calculus rules in convex analysis are
derived as consequences.

If the functions $f_{t}$, $t\in T$, are proper convex and lower
semicontinuous; that is, $\left\{  f_{t},\text{ }t\in T\right\}  \subset
\Gamma_{0}(X),$ and we additionally assume that the relative interior of the
effective domain of $f$ is non-empty, i.e. $\operatorname*{ri}%
(\operatorname*{dom}f)\neq\emptyset,$ in \cite[Lemma 3]{HLZ08} it is
established that%
\begin{equation}
\partial f(x)=%
{\textstyle\bigcap\limits_{\varepsilon>0}}
\overline{\operatorname*{co}}\left(
{\textstyle\bigcup_{t\in T_{\varepsilon}(x)}}
\partial_{\varepsilon}f_{t}(x)+\mathrm{N}_{\operatorname*{dom}f}(x)\right)
,\quad\text{for all }x\in\operatorname*{dom}f, \label{int-formula}%
\end{equation}
where $\overline{\operatorname*{co}}$ stands for the $w^{\ast}$-closed convex
hull, $\partial_{\varepsilon}f_{t}(x)$ is the $\varepsilon$-subdifferential of
$f_{t}$ at $x$, and
\[
T_{\varepsilon}(x):=\{t\in T:\ f_{t}(x)\geq f(x)-\varepsilon\}.
\]
In \cite[Theorem 4]{LoVo10}, formula (\ref{int-formula}) is also derived under
different assumptions, namely if\textbf{ }$\operatorname*{cone}%
(\operatorname*{dom}f-x))$\textbf{ }is closed or\textbf{ }$\operatorname*{ri}%
(\operatorname*{cone}(\operatorname*{dom}f-x))\neq\emptyset,$\textbf{
}where\textbf{ }$\operatorname*{cone}(A)$\textbf{ }is the convex cone
generated by\textbf{ }$A$\textbf{.}

When these interiority/closedness assumptions are removed, the price that has
to be paid is the need of involving the family%
\[
\mathcal{F}(x):=\{L\subset X:\ L\text{ is a finite-dimensional linear subspace
such that }x\in L\}.
\]
In this very general framework, the following characterization is established
in \cite[Theorem 4]{HLZ08}:
\begin{equation}
\partial f(x)=%
{\textstyle\bigcap\limits_{L\in\mathcal{F}(x),\varepsilon>0}}
\overline{\operatorname*{co}}\left(
{\textstyle\bigcup_{t\in T_{\varepsilon}(x)}}
\partial_{\varepsilon}f_{t}(x)+\mathrm{N}_{L\cap\operatorname*{dom}%
f}(x)\right)  ,\ \text{for all }x\in\operatorname*{dom}f. \label{int-formulab}%
\end{equation}
Observe that $\operatorname*{ri}(L\cap\operatorname*{dom}f)\neq\emptyset$. The
reader will find related formulas\ in \cite{LiNg11}.

In the so-called compact setting the following result, involving only the
active functions at the reference point,\textbf{ }is established in
\cite[Theorem 3.8]{CHL19b} under the standard hypothesis (\ref{SH}); i.e., $T$
is compact and the mappings $t\mapsto f_{t}(z),$ $z\in X,$ are upper
semicontinuous (usc, in brief):%
\[
\partial f(x)=%
{\textstyle\bigcap\nolimits_{L\in\mathcal{F}(x),\text{ }\varepsilon>0}}
\overline{\operatorname*{co}}\left(
{\textstyle\bigcup\nolimits_{t\in T(x)}}
\partial_{\varepsilon}f_{t}(x)+\mathrm{N}_{L\cap\operatorname*{dom}%
f}(x)\right)  .
\]

One way to get rid of these normal cones is to impose additional assumptions
as the finiteness and continuity of $f$ at $x$, in which case
(\ref{int-formulab}) gives rise to (\cite[Corollary 10]{HLZ08}; see, also,
\cite{Vo93}, for normed spaces):
\[
\partial f(x)=%
{\textstyle\bigcap\nolimits_{\varepsilon>0}}
\overline{\operatorname*{co}}(%
{\textstyle\bigcup\nolimits_{t\in T_{\varepsilon}(x)}}
\partial_{\varepsilon}f_{t}(x)).
\]

Lemma 5 in \cite{HLZ08} yields some characterizations of $\mathrm{N}%
_{\operatorname*{dom}f}(x).$ Precisely,
\begin{align}
x^{\ast}\in\mathrm{N}_{\operatorname*{dom}f}(x)  &  \Leftrightarrow(x^{\ast
},\left\langle x^{\ast},x\right\rangle )\in\left[  \overline
{\operatorname*{co}}\left(  \cup_{t\in T}\operatorname*{gph}f_{t}^{\ast
}\right)  \right]  _{\infty}\nonumber\\
&  \Leftrightarrow(x^{\ast},\left\langle x^{\ast},x\right\rangle )\in\left[
\overline{\operatorname*{co}}\left(  \cup_{t\in T}\operatorname*{epi}%
f_{t}^{\ast}\right)  \right]  _{\infty}, \label{n2}%
\end{align}
where $\operatorname*{gph}f_{t}^{\ast}$ and $\operatorname*{epi}f_{t}^{\ast}$
represent the graph and the epigraph of the conjugate of $f_{t},$
respectively, and $\left[  \cdot\right]  _{\infty}$ defines\ the recession
cone. In the linear case, i.e. if $f(x):=\sup\{\langle a_{t}^{\ast}%
,x\rangle-b_{t}:\ t\in T\},$ with $a_{t}^{\ast}\in X^{\ast}$ and $b_{t}%
\in\mathbb{R}$, we get (\cite[Corollary 7]{HLZ08})%
\[
x^{\ast}\in\mathrm{N}_{\operatorname*{dom}f}(x)\Leftrightarrow(x^{\ast
},\left\langle x^{\ast},x\right\rangle )\in\left[  \overline
{\operatorname*{co}}\left(  (\theta,0)\cup\{(a_{t}^{\ast},b_{t}),t\in
T\}\right)  \right]  _{\infty},
\]
where $\theta$ is the origin in $X^{\ast}.$

Another interesting problem in optimization consists of characterizing the
normal cone to sublevel sets (see, e.g. \cite{HaSv17, CaTh14} and references
therein). Observe, for instance, that if $g\in\Gamma_{0}(X)$ and $\left[
g\leq0\right]  $ is the $0$-sublevel set, by taking $f:=\sup\nolimits_{\alpha
\geq0}(\alpha g)=\mathrm{I}_{\left[  g\leq0\right]  },$ we obtain that%
\[
\mathrm{N}_{\left[  g\leq0\right]  }(x)=\mathrm{N}_{\operatorname*{dom}f}(x).
\]

The main contribution of this paper consists of formulating alternative
characterizations of $\partial f(x),$ relying exclusively on the data
functions and not on any normal cone. In other words, the normal cone
$\mathrm{N}_{L\cap\operatorname*{dom}f}(x)$ does not appear explicitly in the
new subdifferential\ formulas, and consequently, there is no need of
intersecting over finite-dimensional subspaces $L$, as in previous quoted
works. Extensions to the non-compact framework will be investigated in a
forthcoming work, using different approaches including\ well-known
qualifications, like the strong CHIP, SECQ, linear regularity,
Farkas-Minkowski, etc. (\cite{DGL06, DGLS07, HW10, LiNg05, LiNgP07}).

The structure of the paper is the following. After Section \ref{sec2} devoted
to notation and preliminary results, in Section \ref{sec3} new
characterizations of $\mathrm{N}_{\operatorname*{dom}f}(x)$ are given in terms
exclusively of $\partial_{\varepsilon}f_{t}(x),$ $t\in T,$ which are
independent of $\varepsilon$ and much simpler than those in (\ref{n2}). The
main result in this section is Theorem \ref{p1}.\textbf{ }Based on the results
established in Section \ref{sec3}, Theorems \ref{t1} and \ref{t1bis}\textbf{
}in Section \ref{secsub} provide new formulas for the subdifferential of the
supremum, $\partial f(x),$ involving both, the\textbf{ }active functions at
the reference point $x$, and also the rest of the functions but affected by a
multiplying\textbf{ }parameter. Finally, new optimality conditions for the
convex optimization problem with infinitely many constraints are proposed.

\section{Notation and preliminary results\label{sec2}}

Let\ $X$ be a (real) separated locally convex space (lcs, for short), whose
topological dual space, $X^{\ast},$ is endowed with the $w^{\ast}$-topology;
hence, $X^{\ast\ast}:=(X^{\ast})^{\ast}\equiv X.\ $The spaces $X$ and
$X^{\ast}$ are paired in duality by the bilinear form $(x^{\ast},x)\in
X^{\ast}\times X\mapsto\langle x^{\ast},x\rangle:=x^{\ast}(x).$ The zero
vectors in $X$ and $X^{\ast}$ are denoted by $\theta.$ The basis of closed,
convex and balanced neighborhoods of\textbf{ }$\theta,$ in both $X$ and
$X^{\ast},$ called\textbf{ }$\theta$\textbf{-}neighborhoods\textbf{,} is
represented by\textbf{ }$\mathcal{N}$. We use the notation $\overline
{\mathbb{R}}:=\mathbb{R}\cup\{-\infty,+\infty\}$ and $\mathbb{R}_{\infty
}:=\mathbb{R}\cup\{+\infty\}$, and adopt the conventions\emph{ }$\left(
+\infty\right)  +(-\infty)=\left(  -\infty\right)  +(+\infty)=+\infty,$
$0(+\infty):=+\infty.$

Given $k\geq1,$ we denote
\begin{align*}
\Delta_{k}  &  :=\left\{  (\lambda_{1},\cdots,\lambda_{k})\geq0:\lambda
_{1}+\cdots+\lambda_{k}=1\right\}  ,\\
\Delta_{k}^{+}  &  :=\left\{  (\lambda_{1},\cdots,\lambda_{k})\in\Delta
_{k}:\lambda_{i}>0,\text{ }i=1,\cdots,k\right\}  .
\end{align*}

Given two sets $A$ and $B$ in $X$ (or in $X^{\ast}$), we define the Minkowski
sum by
\begin{equation}
A+B:=\{a+b:\ a\in A,\text{ }b\in B\},\quad A+\emptyset=\emptyset+A=\emptyset,
\label{mincov}%
\end{equation}
and, if\ $\Lambda\subset\mathbb{R},$
\[
\Lambda A:=\left\{  \lambda a:\text{ }\lambda\in\Lambda,\text{ }a\in
A\right\}  ,\quad\Lambda\emptyset=\emptyset A=\emptyset,
\]
in particular, we write $\lambda A:=\left\{  \lambda\right\}  A,$ $\lambda
\in\mathbb{R}.$

By $\operatorname*{co}(A)$ and $\operatorname*{cone}(A)$, we denote the
\emph{convex} and the \emph{conical convex} \emph{hulls }of the nonempty set
$A$, respectively. In the topological side, $\operatorname*{cl}(A)$ and
$\overline{A}$ are indistinctly used for denoting the \emph{closure }of $A$.
When\textbf{ }$A\subset X^{\ast}$, the closure is taken with respect to
the\textbf{ }$w^{\ast}$-topology, unless something else is explicitly stated.

Associated with a nonempty set $A\subset X,$ we define\ the\textbf{
}\emph{negative dual cone} and the \emph{orthogonal} \emph{subspace} of $A$ as
follows
\begin{align*}
A^{-}  &  :=\left\{  x^{\ast}\in X^{\ast}:\ \langle x^{\ast},x\rangle
\leq0\text{ for all }x\in A\right\}  ,\\
A^{\perp}  &  :=(-A^{-})\cap A^{-}=\left\{  x^{\ast}\in X^{\ast}:\ \langle
x^{\ast},x\rangle=0\text{ for all }x\in A\right\}  ,
\end{align*}
respectively. Observe that $A^{-}=(\overline{\operatorname*{cone}}(A))^{-}.$
These concepts are defined similarly for sets in $X^{\ast}.$ The so-called
\emph{bipolar theorem} establishes that
\begin{equation}
A^{--}:=(A^{-})^{-}=\overline{\operatorname*{cone}}(A). \label{bipolar}%
\end{equation}

If $A\subset X,$ we define the \emph{normal cone}\ to $A$ at $x$ by
\[
\mathrm{N}_{A}(x):=\left\{
\begin{array}
[c]{ll}%
(A-x)^{-}, & \text{if \ }x\in A,\\
\emptyset, & \text{if }x\in X\setminus A.
\end{array}
\right.
\]

If\textbf{ }$A\neq\emptyset$ is convex and closed\textbf{,} $A_{\infty}$
represents its \emph{recession cone} defined by
\[
A_{\infty}:=\left\{  y\in X:\ x+\lambda y\in A\text{ for some }x\in A\text{
and all }\lambda\geq0\right\}  .
\]

Given a function $f:X\longrightarrow\overline{\mathbb{R}}$, its
\emph{(effective)} \emph{domain} and \emph{epigraph }are, respectively,
\[
\operatorname*{dom}f:=\{x\in X:\ f(x)<+\infty\},
\]
and%
\[
\operatorname*{epi}f:=\{(x,\lambda)\in X\times\mathbb{R}:\ f(x)\leq\lambda\}.
\]
We say that $f$ is \emph{proper} when\ $\operatorname*{dom}f\neq\emptyset$ and
$f(x)>-\infty$ for all $x\in X$. The \emph{closed hull} of $f$ is the function
$\operatorname*{cl}f:X\longrightarrow\overline{\mathbb{R}}$ whose epigraph is
$\operatorname*{cl}\left(  \operatorname*{epi}f\right)  .$ Moreover,
\begin{equation}
\left(  \operatorname*{cl}f\right)  (x)=\liminf_{x^{\prime}\rightarrow
x}f(x^{\prime})=\sup_{V\in\mathcal{N}}\inf\left\{  f(x^{\prime}):\ x^{\prime
}\in x+V\right\}  . \label{cl}%
\end{equation}

The \emph{convex hull} of $f,$ $\operatorname*{co}f:X\longrightarrow
\overline{\mathbb{R}},$ is the largest convex function which is dominated by
$f.$ Equivalently\textbf{, }%
\begin{align}
(\operatorname*{co}f)(x)  &  =\inf\left\{  \mu:\ (x,\mu)\in\operatorname*{co}%
\left(  \operatorname*{epi}f\right)  \right\} \label{m7}\\
&  =\inf\left\{  \sum_{i=1}^{k}\lambda_{i}f(x_{i}):\ \lambda\in\Delta_{k}%
^{+},\ \sum_{i=1}^{k}\lambda_{i}x_{i}=x,\text{ }k\in\mathbb{N}\right\}  .
\label{m8}%
\end{align}
The \emph{closed convex hull} of $f$ is the convex lower semicontinuous (lsc,
in brief) function $\overline{\operatorname*{co}}f:X\longrightarrow
\overline{\mathbb{R}}$ such that%
\[
\operatorname*{epi}(\overline{\operatorname*{co}}f)=\overline
{\operatorname*{co}}(\operatorname*{epi}f).
\]
Obviously, $\overline{\operatorname*{co}}f\leq\operatorname*{cl}f\leq f.$

Given\ $x\in X$ and $\varepsilon\in\mathbb{R},$ the $\varepsilon
$-\emph{subdifferential }(or the \emph{approximate subdifferential})\emph{ }of
$f$ at $x$ is
\begin{equation}
\partial_{\varepsilon}f(x)=\{x^{\ast}\in X^{\ast}:\ f(y)\geq f(x)+\langle
x^{\ast},y-x\rangle-\varepsilon\text{ \ for all }y\in X\}, \label{defsub}%
\end{equation}
when $x\in\operatorname*{dom}f,$ and $\partial_{\varepsilon}f(x):=\emptyset$
when $f(x)\notin\mathbb{R}$ or $\varepsilon<0.$ The \emph{subdifferential} of
$f$ at $x$ is $\partial f(x):=\partial_{0}f(x)$.

We shall use the following relation (e.g., \cite[Exercise 2.23]{Za02})%
\begin{equation}
\mathrm{N}_{\operatorname*{dom}f}(x)=\left(  \partial_{\varepsilon
}f(x)\right)  _{\infty}, \label{marco29}%
\end{equation}
where\ $f\in\Gamma_{0}(X),$ $x\in\operatorname*{dom}f,$ and $\varepsilon>0.$

The \emph{Fenchel conjugate} of $f$ is the function $f^{\ast}:X^{\ast
}\longrightarrow\overline{\mathbb{R}}\mathbb{\ }$given by
\[
f^{\ast}(x^{\ast}):=\sup\{\left\langle x^{\ast},x\right\rangle -f(x):\ x\in
X\}.
\]

It is well-known that
\begin{equation}
\partial_{\varepsilon}f(x)=\{x^{\ast}\in X^{\ast}:\ f(x)+f^{\ast}(x^{\ast
})\leq\left\langle x^{\ast},x\right\rangle +\varepsilon\}, \label{marco31}%
\end{equation}
for all $\varepsilon\geq0,$ and
\begin{equation}
\partial f(x)=\cap_{\varepsilon>0}\partial_{\varepsilon}f(x). \label{vo}%
\end{equation}
The following lemma gives a slight extension of the last relation, which is
used later on.

\begin{lem}
\label{lemvo}Consider a function $f\in\Gamma_{0}(X)$ and suppose that $x$ is a
minimizer of $f.\ $Then, for every\ $M\geq0,$%
\begin{equation}
\partial f(x)=\cap_{\varepsilon>0}\overline{\operatorname*{co}}\left(
\partial_{\varepsilon}f(x)\cup\varepsilon\partial_{\varepsilon+M}f(x)\right)
=\cap_{\varepsilon>0}\overline{\operatorname*{co}}\left(  \partial
_{\varepsilon}f(x)\cup\varepsilon\partial_{\varepsilon+M}f^{+}(x)\right)  ,
\label{vob}%
\end{equation}
where $f^{+}:=\max\left\{  f,0\right\}  $ is the positive part of $f.$
\end{lem}

\begin{dem}
The inclusions
\begin{equation}
\partial f(x)\subset\cap_{\varepsilon>0}\overline{\operatorname*{co}}\left(
\partial_{\varepsilon}f(x)\cup\varepsilon\partial_{\varepsilon+M}f(x)\right)
\label{l1}%
\end{equation}
and
\begin{equation}
\partial f(x)\subset\cap_{\varepsilon>0}\overline{\operatorname*{co}}\left(
\partial_{\varepsilon}f(x)\cup\varepsilon\partial_{\varepsilon+M}%
f^{+}(x)\right)  \label{l2}%
\end{equation}
follow easily from (\ref{vo}). We only need to prove the opposite inclusion in
(\ref{l1}), since that the same argument is valid for such an inclusion
in\ (\ref{l2}). Take $x^{\ast}$ in the right-hand side of (\ref{l1}). Then,
for each fixed $\varepsilon>0,$
\[
x^{\ast}=\lim_{i}(\lambda_{i}y_{i}^{\ast}+(1-\lambda_{i})\varepsilon
z_{i}^{\ast}),
\]
for some nets $(\lambda_{i})_{i}\subset\left[  0,1\right]  ,$ $(y_{i}^{\ast
})_{i}\subset\partial_{\varepsilon}f(x)$ and $(z_{i}^{\ast})_{i}%
\subset\partial_{\varepsilon+M}f(x).$ Thus, for each\ $y\in\operatorname*{dom}%
f,$
\begin{align*}
\left\langle x^{\ast},y-x\right\rangle  &  =\lim_{i}\left\langle \lambda
_{i}y_{i}^{\ast}+(1-\lambda_{i})\varepsilon z_{i}^{\ast},y-x\right\rangle \\
&  \leq\limsup_{i}\left(  \lambda_{i}(f(y)-f(x)+\varepsilon)+(1-\lambda
_{i})\varepsilon(f(y)-f(x)+\varepsilon+M)\right) \\
&  \leq f(y)-f(x)+\varepsilon+\varepsilon(f(y)-f(x)+\varepsilon+M),
\end{align*}
as $f(y)\geq f(x).$ Finally, since the last inequality holds for all
$\varepsilon>0,$ when $\varepsilon\downarrow0$ we obtain
\[
\left\langle x^{\ast},y-x\right\rangle \leq f(y)-f(x)\text{ for all }y\in X,
\]
which shows that $x^{\ast}\in\partial f(x).$
\end{dem}

The \emph{support} and the \emph{indicator} functions of $A\subset X$ are,
respectively,%
\[
\sigma_{A}(x^{\ast}):=\sup\{\langle x^{\ast},x\rangle:\ x\in A\},\text{
}x^{\ast}\in X^{\ast},
\]
with $\sigma_{\emptyset}\equiv-\infty$, and
\[
\mathrm{I}_{A}(x):=\left\{
\begin{array}
[c]{ll}%
0 & \text{if }x\in A,\\
+\infty & \text{if }x\in X\setminus A.
\end{array}
\right.
\]

It is known that, if\textbf{ }$A$ is a closed convex set\textbf{,}
\begin{equation}
A_{\infty}=\left(  \operatorname*{dom}\sigma_{A}\right)  ^{-}, \label{marco17}%
\end{equation}
or equivalently, by using (\ref{bipolar}),%
\begin{equation}
\left(  A_{\infty}\right)  ^{-}=\operatorname*{cl}(\operatorname*{dom}%
\sigma_{A}). \label{domsup}%
\end{equation}

\begin{lem}
\label{consum-1} Consider nonempty sets $A$ and $A_{1},$ $\cdots,A_{k}$ in
$X,$ $k\geq2.$ Then
\begin{equation}
\left[  \overline{\operatorname*{co}}(A\cup\left(  \cup_{i=1,\cdots,k}%
A_{k}\right)  )\right]  _{\infty}=\left(  \overline{\operatorname*{co}}%
(A\cup(A_{1}+\cdots+A_{k}))\right)  _{\infty}. \label{marco13}%
\end{equation}

\end{lem}

\begin{dem}
Obviously,
\begin{align*}
\operatorname*{dom}\mathrm{\sigma}_{A\cup\left(  \cup_{i=1,\cdots,k}%
A_{k}\right)  }  &  =\operatorname*{dom}\left(  \max\left\{  \mathrm{\sigma
}_{A};\text{ }\mathrm{\sigma}_{A_{i}},i=1,\cdots,k\right\}  \right) \\
&  =\operatorname*{dom}\left(  \max\left\{  \mathrm{\sigma}_{A},\text{
}\mathrm{\sigma}_{A_{1}}+\cdots+\mathrm{\sigma}_{A_{k}}\right\}  \right) \\
&  =\operatorname*{dom}\left(  \max\left\{  \mathrm{\sigma}_{A},\text{
}\mathrm{\sigma}_{A_{1}+\cdots+A_{k}}\right\}  \right) \\
&  =\operatorname*{dom}\left(  \mathrm{\sigma}_{A\cup(A_{1}+\cdots+A_{k}%
)}\right)  ,
\end{align*}
and, by (\ref{domsup}),
\begin{align*}
\left(  \left[  \overline{\operatorname*{co}}(A\cup\left(  \cup_{i=1,\cdots
,k}A_{k}\right)  )\right]  _{\infty}\right)  ^{-}  &  =\operatorname*{cl}%
(\operatorname*{dom}\mathrm{\sigma}_{\overline{\operatorname*{co}}(A\cup
(\cup_{i=1,\cdots,k}A_{k}))})=\operatorname*{cl}(\operatorname*{dom}%
\mathrm{\sigma}_{A\cup(\cup_{i=1,\cdots,k}A_{k})})\\
&  =\operatorname*{cl}(\operatorname*{dom}\left(  \mathrm{\sigma}_{A\cup
(A_{1}+\cdots+A_{k})}\right)  )=\left[  (\overline{\operatorname*{co}}%
(A\cup(A_{1}+\cdots+A_{k})))_{\infty}\right]  ^{-}.
\end{align*}
Thus,\ (\ref{marco13}) follows from (\ref{bipolar}).
\end{dem}

\begin{lem}
\label{lemconsum0} Consider a family of nonempty sets $\left\{  A_{t},\ t\in
T_{1}\cup T_{2}\right\}  \subset X,$\ where $T_{1}$ and $T_{2}$ are disjoint
nonempty sets\textbf{. }Then for every $m>0$ we have
\begin{align}
\left[  \overline{\operatorname*{co}}\left(  \cup_{t\in T_{1}\cup T_{2}}%
A_{t}\right)  \right]  _{\infty}  &  =\left[  \overline{\operatorname*{co}%
}\left(  (\cup_{t\in T_{1}}A_{t})\cup(\cup_{t\in T_{2}}mA_{t})\right)
\right]  _{\infty}\nonumber\\
&  =\left[  \overline{\operatorname*{co}}\left(  \cup_{t_{1}\in T_{1},t_{2}\in
T_{2}}\left(  A_{t_{1}}+mA_{t_{2}}\right)  \right)  \right]  _{\infty}.
\label{23}%
\end{align}

\end{lem}

\begin{dem}
Denote $T:=T_{1}\cup T_{2}$ and $A:=\overline{\operatorname*{co}}\left(
\cup_{t\in T}A_{t}\right)  .$ Then the functions $\varphi_{1}:=\sup_{t\in
T_{1}}\mathrm{\sigma}_{A_{t}},$ $\varphi_{2}:=\sup_{t\in T_{2}}\mathrm{\sigma
}_{mA_{t}}$ satisfy
\begin{align}
\varphi_{1}+\varphi_{2}  &  =\sup_{t_{1}\in T_{1},t_{2}\in T_{2}%
}(\mathrm{\sigma}_{A_{t_{1}}}+\mathrm{\sigma}_{mA_{t_{2}}})\nonumber\\
&  =\sup\limits_{t_{1}\in T_{1},t_{2}\in T_{2}}\mathrm{\sigma}_{A_{t_{1}%
}+mA_{t_{2}}}=\mathrm{\sigma}_{\cup_{t_{1}\in T_{1},t_{2}\in T_{2}}(A_{t_{1}%
}+mA_{t_{2}})}, \label{i0}%
\end{align}
and similarly%
\[
\max\left\{  \varphi_{1},\varphi_{2}\right\}  =\mathrm{\sigma}_{(\cup_{t\in
T_{1}}A_{t})\cup(\cup_{t\in T_{2}}mA_{t})}.
\]
Since
\[
\operatorname*{dom}(\varphi_{1}+\varphi_{2})=\operatorname*{dom}(\max\left\{
\varphi_{1},\varphi_{2}\right\}  )=\operatorname*{dom}(\max\left\{
\varphi_{1},m^{-1}\varphi_{2}\right\}  ),
\]
(\ref{domsup}) yields
\begin{align*}
\left(  \left[  \overline{\operatorname*{co}}\left(  \cup_{t_{1}\in
T_{1},\text{ }t_{2}\in T_{2}}\left(  A_{t_{1}}+mA_{t_{2}}\right)  \right)
\right]  _{\infty}\right)  ^{-}  &  =\left(  \left[  \overline
{\operatorname*{co}}\left(  (\cup_{t\in T_{1}}A_{t})\cup(\cup_{t\in T_{2}%
}mA_{t})\right)  \right]  _{\infty}\right)  ^{-}\\
&  =\left(  \left[  \overline{\operatorname*{co}}\left(  (\cup_{t\in T_{1}%
}A_{t})\cup(\cup_{t\in T_{2}}A_{t})\right)  \right]  _{\infty}\right)  ^{-},
\end{align*}
and we are done thanks to\ (\ref{bipolar}).
\end{dem}

The following lemma provides\ the $\varepsilon$-subdifferential of the
positive part of convex functions. It can be derived\ from \cite[Corollary
2.8.11]{Za02} but we prefer to give here a simple alternative proof based on
\cite[Lemma 1]{LoVo11}.

\begin{lem}
\label{hamu3l} Consider a function\ $f\in\Gamma_{0}(X)$ and let $x\in
\operatorname*{dom}f.$ Then, for every $\varepsilon\geq0,$\ we have
\begin{equation}
\partial_{\varepsilon}f^{+}(x)=%
{\textstyle\bigcup\nolimits_{0\leq\lambda\leq1}}
\partial_{\varepsilon+\lambda f(x)-f^{+}(x)}(\lambda f)(x). \label{hamu3b}%
\end{equation}
where $(\lambda f)(z):=\lambda f(z),$ $z\in X$ (with $0f=\mathrm{I}%
_{\operatorname*{dom}f}).$
\end{lem}

\begin{dem}
By \cite[Lemma 1]{LoVo11} we have that
\[
(f^{+})^{\ast}=\operatorname*{cl}\left(  \inf\nolimits_{\lambda\in\left[
0,1\right]  }(\lambda f+\mathrm{I}_{\operatorname*{dom}f})^{\ast}\right)
\equiv\operatorname*{cl}\left(  \inf\nolimits_{\lambda\in\left[  0,1\right]
}(\lambda f)^{\ast}\right)  .
\]
We take $x^{\ast}\in\partial_{\varepsilon}f^{+}(x)$ and fix\ $n\geq1.$ Then,
by (\ref{marco31}),
\begin{align*}
\operatorname*{cl}\left(  \inf\nolimits_{\lambda\in\left[  0,1\right]
}(\lambda f)^{\ast}\right)  (x^{\ast})+f^{+}(x)  &  =(f^{+})^{\ast}(x^{\ast
})+f^{+}(x)\\
&  \leq\left\langle x^{\ast},x\right\rangle +\varepsilon<\left\langle x^{\ast
},x\right\rangle +\varepsilon+\frac{1}{n},
\end{align*}
and so, there exists a net $(x_{n,i}^{\ast})_{i}$ which ($w^{\ast}$-)converges
to $x^{\ast}$ and satisfies
\begin{align*}
\lim_{i}\left(  \inf\nolimits_{\lambda\in\left[  0,1\right]  }(\lambda
f)^{\ast}\right)  (x_{n,i}^{\ast})+f^{+}(x)  &  <\left\langle x^{\ast
},x\right\rangle +\varepsilon+\frac{1}{n}\\
&  =\lim_{i}\left\langle x_{n,i}^{\ast},x\right\rangle +\varepsilon+\frac
{1}{n}.
\end{align*}
Without loss of generality, we may suppose that
\[
\inf\nolimits_{\lambda\in\left[  0,1\right]  }(\lambda f)^{\ast}(x_{n,i}%
^{\ast})+f^{+}(x)<\left\langle x_{n,i}^{\ast},x\right\rangle +\varepsilon
+\frac{1}{n}\text{ for all }i.
\]
Hence, there exists some $\lambda_{n,i}\in\left[  0,1\right]  $ such that
\[
(\lambda_{n,i}f)^{\ast}(x_{n,i}^{\ast})+f^{+}(x)<\left\langle x_{n,i}^{\ast
},x\right\rangle +\varepsilon+\frac{1}{n}\text{ for all }i,
\]
and consequently, using a diagonal argument, we find $(x_{n,i(n)}^{\ast})_{n}$
and $(\lambda_{n,i(n)})_{n}$ such that $(x_{n,i(n)}^{\ast})_{n}$ ($w^{\ast}%
$-)converges\ to $x^{\ast}$ and
\[
(\lambda_{n,i(n)}f)^{\ast}(x_{n,i(n)}^{\ast})+f^{+}(x)<\left\langle
x_{n,i(n)}^{\ast},x\right\rangle +\varepsilon+\frac{1}{n}\text{ \ for all
}n\geq1.
\]
Equivalently, we have that
\[
(\lambda_{n,i(n)}f)^{\ast}(x_{n,i(n)}^{\ast})+(\lambda_{n,i(n)}%
f)(x)<\left\langle x_{n,i(n)}^{\ast},x\right\rangle +(\lambda_{n,i(n)}%
f)(x)-f^{+}(x)+\varepsilon+\frac{1}{n},
\]
so that $(\lambda_{n,i(n)}f)(x)-f^{+}(x)+\varepsilon+\frac{1}{n}>0,$ because
$\left\langle x_{n,i(n)}^{\ast},x\right\rangle \leq(\lambda_{n,i(n)}f)^{\ast
}(x_{n,i(n)}^{\ast})+(\lambda_{n,i(n)}f)(x),$ and\ $x_{n,i(n)}^{\ast}%
\in\partial_{(\lambda_{n,i(n)}f)(x)-f^{+}(x)+\varepsilon+\frac{1}{n}}%
(\lambda_{n,i(n)}f)(x).$ Therefore, for each $y\in X,$\textbf{ }%
\begin{align*}
\left\langle x_{n,i(n)}^{\ast},y-x\right\rangle  &  \leq(\lambda
_{n,i(n)}f)(y)-(\lambda_{n,i(n)}f)(x)+((\lambda_{n,i(n)}f)(x)-f^{+}%
(x)+\varepsilon+\frac{1}{n})\\
&  =(\lambda_{n,i(n)}f)(y)-f^{+}(x)+\varepsilon+\frac{1}{n}.
\end{align*}
Since\textbf{ }$(\lambda_{n,i(n)})_{n},$ or a subnet of it\textbf{, }converges
to some $\lambda\in\left[  0,1\right]  $, one gets by taking\ limits on the
inequality, for all\ $y\in X,$
\[
\left\langle x^{\ast},y-x\right\rangle \leq(\lambda f)(y)-f^{+}(x)+\varepsilon
=(\lambda f)(y)-(\lambda f)(x)+((\lambda f)(x)-f^{+}(x)+\varepsilon).
\]
In other words, $x^{\ast}\in%
{\textstyle\bigcup\nolimits_{0\leq\lambda\leq1}}
\partial_{\varepsilon+(\lambda f)(x)-f^{+}(x)}(\lambda f)(x)$ showing that the
inclusion \textquotedblleft$\subset$\textquotedblright\ holds.

To prove the converse inclusion \textquotedblleft$\supset$\textquotedblright%
\ pick $x^{\ast}\in\partial_{\varepsilon+\lambda f(x)-f^{+}(x)}(\lambda f)(x)$
for some $0\leq\lambda\leq1.$ Then, for every $y\in X,$%
\begin{align*}
\left\langle x^{\ast},y-x\right\rangle  &  \leq(\lambda f)(y)-(\lambda
f)(x)+\varepsilon+\lambda f(x)-f^{+}(x)\\
&  =(\lambda f)(y)-f^{+}(x)+\varepsilon\\
&  \leq(\lambda f^{+})(y)-f^{+}(x)+\varepsilon\\
&  \leq f^{+}(y)-f^{+}(x)+\varepsilon,
\end{align*}
and $x^{\ast}\in\partial_{\varepsilon}f^{+}(x).$
\end{dem}

\section{Normal cone to the domain\label{sec3}}

This section is devoted to give a representation\ of the normal cone to the
effective domain of a supremum function by means of the $\varepsilon
$-subdifferential of the data functions. This result will be a key tool\ to
derive new formulas for the subdifferential of the supremum function in
Section \ref{secsub}.

We consider a nonempty family $\left\{  f_{t},\ t\in T\right\}  \subset
\Gamma_{0}(X),$ where $X$ is a given lcs space, and the associated\textbf{
}supremum function
\[
f=\sup_{t\in T}f_{t}.
\]
Our analysis is carried out in\textbf{ }the following standard
framework\textbf{:}%
\begin{equation}
T\text{ is Hausdorff compact and the mappings }t\mapsto f_{t}(z),\text{ }z\in
X,\text{ are usc on }T. \label{SH}%
\end{equation}
Given $x\in\operatorname*{dom}f$ and $\varepsilon\geq0,$ remember that the
$\varepsilon$-active set at $x$ is\
\[
T_{\varepsilon}(x):=\left\{  t\in T:\text{ }f_{t}(x)\geq f(x)-\varepsilon
\right\}  ,\text{ \ }T(x):=T_{0}(x).
\]

We shall need the following lemma, which is valid for any family of convex
functions, not necessarily lsc. We refer to \cite[Page 854]{HL08} for a
finite-dimensional version of this result (without proof).

\begin{lem}
Provided assumption \emph{(\ref{SH})} holds, we have that
\begin{equation}
\operatorname*{dom}f=\cap_{t\in T}\operatorname*{dom}f_{t}, \label{egdom}%
\end{equation}
and, for every $x\in\operatorname*{dom}f,$%
\begin{equation}
\mathbb{R}_{+}(\operatorname*{dom}f-x)=\cap_{t\in T}\mathbb{R}_{+}%
(\operatorname*{dom}f_{t}-x). \label{egconos}%
\end{equation}

\end{lem}

\begin{dem}
Take $z\in\cap_{t\in T}\operatorname*{dom}f_{t}.$ For each $t\in T$ the upper
semicontinuity assumption yields some $m_{t}\geq0$ and a neighborhood $V_{t}$
of $t$ such that
\[
f_{s}(z)\leq m_{t}\text{ \ for all }s\in V_{t}.
\]
Consider $V_{t_{1}},$ $\cdots,V_{t_{k}}$ a finite covering of $T.$ Then, for
each $t\in T,$
\[
f_{t}(z)\leq\max\left\{  m_{i},\text{ }i=1,\cdots,k\right\}  <+\infty,
\]
and so $z\in\operatorname*{dom}f.$

To\ prove the second statement we take $z\in\cap_{t\in T}\mathbb{R}%
_{+}(\operatorname*{dom}f_{t}-x).$ If $z=\theta,$ then we are obviously done.
Otherwise, for each $t\in T$ there exist $\alpha_{t},m_{t}>0$ and $z_{t}%
\in\operatorname*{dom}f_{t}$ such that $z=\alpha_{t}(z_{t}-x)$ and
\[
f_{t}(z_{t})<m_{t}.
\]
By arguing as above, the upper semicontinuity assumption yields some
neighborhood $V_{t}$ of $t$ such that
\[
f_{s}(z_{t})<m_{t}\text{ \ for all }s\in V_{t}.
\]
Consider $V_{t_{1}},$ $\cdots,V_{t_{k}},$ a finite covering of $T,$ so
that\ $z_{t_{i}}=\alpha_{t_{i}}^{-1}z+x\in\operatorname*{dom}f_{s}$ for all
$s\in V_{t_{i}}.$ Next, for\ $\bar{\alpha}:=\max\left\{  \alpha_{t_{i}},\text{
}i=1,\cdots,k\right\}  $ we obtain that
\[
\bar{\alpha}^{-1}z=\bar{\alpha}^{-1}(\alpha_{t_{i}}(z_{t_{i}}-x))\in
\bar{\alpha}^{-1}\alpha_{t_{i}}(\cap_{s\in V_{t_{i}}}(\operatorname*{dom}%
f_{s}-x))\subset\cap_{s\in V_{t_{i}}}(\operatorname*{dom}f_{s}-x),\text{
\ }i=1,\cdots,k,
\]
where the last inclusion comes from the convexity of the set $\cap_{s\in
V_{t_{i}}}(\operatorname*{dom}f_{s}-x)$ and the fact that $\theta\in\cap_{s\in
V_{t_{i}}}(\operatorname*{dom}f_{s}-x)$. Hence,%
\begin{align*}
\bar{\alpha}^{-1}z  &  \in\cap_{i=1,\cdots,k}\cap_{s\in V_{t_{i}}%
}(\operatorname*{dom}f_{s}-x)\\
&  =\cap_{s\in V_{t_{i}},i=1,\cdots,k}(\operatorname*{dom}f_{s}-x)=\cap_{t\in
T}(\operatorname*{dom}f_{t}-x),
\end{align*}
and the first statement of the lemma leads us to
\[
z\in\bar{\alpha}(\cap_{t\in T}(\operatorname*{dom}f_{t}-x))\subset
\mathbb{R}_{+}\left(  \cap_{t\in T}(\operatorname*{dom}f_{t}-x)\right)
=\mathbb{R}_{+}(\operatorname*{dom}f-x).
\]

\end{dem}

We give now the main result of this section. When $T$ is a singleton, it
reduces to (\ref{marco29}). For instance, if we apply (\ref{marco29}) to the
supremum function $f$ we obtain that $\mathrm{N}_{\operatorname*{dom}%
f}(x)=\left(  \partial_{\varepsilon}f(x)\right)  _{\infty}.$ Thus, one may
think of using one of the known formulas of the $\varepsilon$-subdifferential
of the supremum function $f$ like in \cite[Theorem 2]{HaSv17} (see, also,
\cite[Theorem 5]{PP}, \cite{PPBis}),\textbf{ }but these formulas involve
approximate subdifferentials $\partial_{\beta_{t}}f_{t}(x)$ with possibly very
large parameters $\beta_{t},$ which are\textbf{ }out of control.

\begin{theo}
\label{p1} Consider $x\in\operatorname*{dom}f$ and let $0<\varepsilon_{t}%
\leq1,$ $t\in T,\ $be such that\
\begin{equation}
\inf_{t\in T}(\varepsilon_{t}f_{t})(x)>-\infty. \label{cb2}%
\end{equation}
Then for every $\varepsilon>0$ we have that
\begin{equation}
\left[  \overline{\operatorname*{co}}\left(
{\textstyle\bigcup\nolimits_{t\in T}}
\partial_{\varepsilon}(\varepsilon_{t}f_{t})(x)\right)  \right]  _{\infty
}\subset\mathrm{N}_{\operatorname*{dom}f}(x), \label{marco37}%
\end{equation}
and, if the standard assumption \emph{(\ref{SH}) }holds, then
\[
\mathrm{N}_{\operatorname*{dom}f}(x)=\left[  \overline{\operatorname*{co}%
}\left(
{\textstyle\bigcup\nolimits_{t\in T}}
\partial_{\varepsilon}(\varepsilon_{t}f_{t})(x)\right)  \right]  _{\infty}.
\]

\end{theo}

\begin{dem}
We fix $\varepsilon>0$ and denote
\[
E_{\varepsilon}:=%
{\textstyle\bigcup\nolimits_{t\in T}}
\partial_{\varepsilon}(\varepsilon_{t}f_{t})(x);
\]
hence, $E_{\varepsilon}\neq\emptyset$ as $\left\{  \varepsilon_{t}f_{t},\text{
}t\in T\right\}  \subset\Gamma_{0}(X).$

To establish the inclusion (\ref{marco37}), we take\ $x^{\ast}\in\left[
\overline{\operatorname*{co}}(E_{\varepsilon})\right]  _{\infty}$ and fix
$x_{0}^{\ast}\in E_{\varepsilon}.$ Then for every $\alpha>0$ we have that
$x_{0}^{\ast}+\alpha x^{\ast}\in\overline{\operatorname*{co}}(E_{\varepsilon
})$ and, so, there are nets $(\lambda_{j,1},\cdots,\lambda_{j,k_{j}})\in
\Delta_{k_{j}}^{+},$ $t_{j,1},\cdots,t_{j,k_{j}}\in T,$ and $x_{j,1}^{\ast}%
\in\partial_{\varepsilon}(\varepsilon_{t_{j,1}}f_{t_{j,1}})(x),\cdots,$
$x_{j,k_{j}}^{\ast}\in\partial_{\varepsilon}(\varepsilon_{t_{j,k_{j}}%
}f_{t_{j,k_{j}}})(x)$ such that
\[
x_{0}^{\ast}+\alpha x^{\ast}=\lim_{j}(\lambda_{j,1}x_{j,1}^{\ast}%
+\cdots+\lambda_{j,k_{j}}x_{j,k_{j}}^{\ast}).
\]
Hence, for every fixed $y\in\operatorname*{dom}f,$\
\begin{align*}
\left\langle x_{0}^{\ast}+\alpha x^{\ast},y-x\right\rangle  &  =\lim
\nolimits_{j}\left\langle \lambda_{j,1}x_{j,1}^{\ast}+\cdots+\lambda_{j,k_{j}%
}x_{j,k_{j}}^{\ast},y-x\right\rangle \\
&  \leq\limsup\nolimits_{j}\left(  \sum_{i=1,\cdots,k_{j}}\lambda
_{j,i}(\varepsilon_{t_{j,i}}f_{t_{j,i}}(y)-\varepsilon_{t_{j,i}}f_{t_{j,i}%
}(x)+\varepsilon)\right) \\
&  \leq\limsup\nolimits_{j}\left(  \sum_{i=1,\cdots,k_{j}}\lambda
_{j,i}(\varepsilon_{t_{j,i}}f^{+}(y)-\varepsilon_{t_{j,i}}f_{t_{j,i}%
}(x)+\varepsilon)\right) \\
&  \leq f^{+}(y)-\inf\left\{  \varepsilon_{t}f_{t}(x),\text{ }t\in T\right\}
+\varepsilon,
\end{align*}
and condition (\ref{cb2}) ensures, by dividing by $\alpha$ and making
$\alpha\uparrow+\infty,$ that\
\[
\left\langle x^{\ast},y-x\right\rangle \leq0,
\]
for all $y\in\operatorname*{dom}f=\operatorname*{dom}f^{+},$ that is,
$x^{\ast}\in\mathrm{N}_{\operatorname*{dom}f}(x),$ as we wanted to prove.

Now we assume that the additional condition\ (\ref{SH}) holds. We have to
prove the inclusion
\begin{equation}
\left(  \left[  \overline{\operatorname*{co}}\left(  E_{\varepsilon}\right)
\right]  _{\infty}\right)  ^{-}\subset(\mathrm{N}_{\operatorname*{dom}%
f}(x))^{-}, \label{mn}%
\end{equation}
or equivalently, according to (\ref{domsup}) and the fact that $(\mathrm{N}%
_{\operatorname*{dom}f}(x))^{-}=\operatorname*{cl}(\mathbb{R}_{+}%
(\operatorname*{dom}f-x)),$
\begin{equation}
\operatorname*{cl}(\operatorname*{dom}\mathrm{\sigma}_{E_{\varepsilon}%
})\subset\operatorname*{cl}(\mathbb{R}_{+}(\operatorname*{dom}f-x)).
\label{mn2}%
\end{equation}
Take
\begin{align*}
z\in\operatorname*{dom}\mathrm{\sigma}_{E_{\varepsilon}}  &
=\operatorname*{dom}\left(  \mathrm{\sigma}_{%
{\textstyle\bigcup\nolimits_{t\in T}}
\partial_{\varepsilon}(\varepsilon_{t}f_{t})(x)}\right) \\
&  =\operatorname*{dom}\left(  \sup_{t\in T}\mathrm{\sigma}_{\partial
_{\varepsilon}(\varepsilon_{t}f_{t})(x)}\right)  =\operatorname*{dom}\left(
\sup_{t\in T}(\varepsilon_{t}f_{t})_{\varepsilon}^{\prime}(x;\cdot)\right)  ,
\end{align*}
where $(\varepsilon_{t}f_{t})_{\varepsilon}^{\prime}(x;\cdot)$ represents the
$\varepsilon$-directional derivative of the function $\varepsilon_{t}f_{t}%
\in\Gamma_{0}(X)$ at $x\ $(see \cite[Theorem 2.4.11]{Za02}). Then, by
\cite[Theorem 2.1.14]{Za02},
\begin{align*}
z  &  \in\cap_{t\in T}\operatorname*{dom}((\varepsilon_{t}f_{t})_{\varepsilon
}^{\prime}(x;\cdot))\\
&  =\cap_{t\in T}\mathbb{R}_{+}(\operatorname*{dom}(\varepsilon_{t}f_{t})-x)\\
&  =\cap_{t\in T}\mathbb{R}_{+}(\operatorname*{dom}f_{t}-x),
\end{align*}
and (\ref{egconos}) gives rise to
\[
z\in\cap_{t\in T}\mathbb{R}_{+}(\operatorname*{dom}f_{t}-x)=\mathbb{R}%
_{+}(\operatorname*{dom}f-x)\subset\operatorname*{cl}(\mathbb{R}%
_{+}(\operatorname*{dom}f-x)).
\]
Hence, (\ref{mn2}) holds\ and (\ref{mn}) follows.
\end{dem}

The following corollary, which is straightforward from Theorem \ref{p1}, gives
a practical example for\ the weighting parameters $\varepsilon_{t},$ $t\in T,$
used in\ the above characterization of $\mathrm{N}_{\operatorname*{dom}f}(x).$

\begin{cor}
\label{cp1} Consider $\varepsilon>0,$ $x\in\operatorname*{dom}f$ and denote
\begin{equation}
\varepsilon_{t}:=\left\{
\begin{array}
[c]{ll}%
\frac{-\varepsilon}{2f_{t}(x)-2f(x)+\varepsilon}\mathbf{,} & \text{if }t\in
T\setminus T_{\varepsilon}(x),\\
1, & \text{if }t\in T_{\varepsilon}(x).
\end{array}
\right.  \label{vet}%
\end{equation}
Then $0<\varepsilon_{t}<1$ for all $t\in T_{\varepsilon}(x)\ $and, provided
that \emph{(\ref{SH}) }holds,
\[
\mathrm{N}_{\operatorname*{dom}f}(x)=\left[  \overline{\operatorname*{co}%
}\left(  \left(
{\textstyle\bigcup\nolimits_{t\in T_{\varepsilon}(x)}}
\partial_{\varepsilon}f_{t}(x)\right)  \cup\left(
{\textstyle\bigcup\nolimits_{t\in T\setminus T_{\varepsilon}(x)}}
\partial_{\varepsilon}(\varepsilon_{t}f_{t})(x)\right)  \right)  \right]
_{\infty}.
\]

\end{cor}

\begin{dem}
We may assume that $f(x)=0.$ Then, for every $t\in T\setminus T_{\varepsilon
}(x),$ we have that $2f_{t}(x)+\varepsilon<-\varepsilon$ and, so,
$0<\varepsilon_{t}<1.$ Also, for\ such\ $t\in T\setminus T_{\varepsilon}(x)$
we have that
\[
\varepsilon_{t}f_{t}(x)=\frac{-\varepsilon f_{t}(x)}{2f_{t}(x)+\varepsilon
}>-\varepsilon,
\]
so that $\inf_{t\in T}(\varepsilon_{t}f_{t})(x)\geq-\varepsilon$ and condition
(\ref{cb2}) follows. Thus, the desired conclusion straightforwardly
comes\ from Theorem \ref{p1}.
\end{dem}

The following result is a simple consequence of Theorem \ref{p1}, giving a
characterization of $\mathrm{N}_{\operatorname*{dom}f}(x)\ $by means of the
original\ functions $f_{t}$'s and not the $(\varepsilon_{t}f_{t})$'s.

\begin{cor}
\label{ccor} Consider $x\in\operatorname*{dom}f$ and assume that\emph{ }%
\begin{equation}
\inf_{t\in T}f_{t}(x)>-\infty. \label{cb}%
\end{equation}
Then, under\ \emph{(\ref{SH}),} for every $\varepsilon>0$ we have
\begin{equation}
\mathrm{N}_{\operatorname*{dom}f}(x)=\left[  \overline{\operatorname*{co}%
}\left(
{\textstyle\bigcup\nolimits_{t\in T}}
\partial_{\varepsilon}f_{t}(x)\right)  \right]  _{\infty}. \label{fg}%
\end{equation}

\end{cor}

\begin{dem}
Take $\varepsilon_{t}=1,$ $t\in T,$ in Theorem \ref{p1}.
\end{dem}

Condition (\ref{cb}) obviously holds when\ $T$ is finite. More generally, we
have the following result.

\begin{cor}
Consider $x\in\operatorname*{dom}f$ and assume that \emph{(\ref{SH}) }holds.
If the mapping $t\mapsto f_{t}(x)$ is also\ lsc, then condition
\emph{(\ref{cb})} fulfills. Consequently, \emph{(\ref{fg})} is satisfied.
\end{cor}

\begin{dem}
Since $\left\{  f_{t},t\in T\right\}  \subset\Gamma_{0}(X)$, for each $t\in T$
we have that
\[
+\infty>f(x)\geq f_{t}(x)>f_{t}(x)-1,
\]
and so there exists some neighborhood $V_{t}$ of $t$ such that
\[
f_{s}(x)>f_{t}(x)-1,\text{ \ for all }s\in V_{t}.
\]
But $T$ is compact, and so $T\subset\cup_{i=1}^{k}V_{t_{i}},$ for some
$\left\{  t_{1},\cdots,t_{k}\right\}  \subset T.$ Hence, for each $t\in T,$ we
have that $t\in V_{t_{i_{0}}}$ for some $i\in\left\{  1,\cdots,k\right\}  ,$
so that
\[
f_{t}(x)>f_{t_{i_{0}}}(x)-1\geq\min_{i\in\left\{  1,\cdots,k\right\}
}f_{t_{i}}(x)-1>-\infty.
\]

\end{dem}

The following corollary shows that we can give different values to\ the
parameter $\varepsilon$ and the formula in Theorem \ref{p1} is still valid. It
is an extension to our\textbf{ }current setting of \cite[Lemma 11(ii)]{HLZ08}
dealing with finitely many functions.

\begin{cor}
\label{lemconsum} Assume that hypothesis \emph{(\ref{SH}) }holds. Given
$x\in\operatorname*{dom}f,$ $(\varepsilon_{t})_{t}\subset\left]  0,1\right]  $
satisfying \emph{(\ref{cb2})}, and $\delta_{t}>0,$ $t\in T,$ such that
$0<\inf_{t\in T}\delta_{t}\leq\sup_{t\in T}\delta_{t}<+\infty,$ we have
\[
\mathrm{N}_{\operatorname*{dom}f}(x)=\left[  \overline{\operatorname*{co}%
}\left(  \cup_{t\in T}\partial_{\delta_{t}}(\varepsilon_{t}f_{t})(x)\right)
\right]  _{\infty}.
\]

\end{cor}

\begin{dem}
We denote $\bar{\delta}:=\inf_{t\in T}\delta_{t}$ and $\hat{\delta}%
:=\sup_{t\in T}\delta_{t}.$ Then
\[
\left[  \overline{\operatorname*{co}}\left(  \cup_{t\in T}\partial
_{\bar{\delta}}(\varepsilon_{t}f_{t})(x)\right)  \right]  _{\infty}%
\subset\left[  \overline{\operatorname*{co}}\left(  \cup_{t\in T_{1}}%
\partial_{\delta_{t}}(\varepsilon_{t}f_{t})(x)\right)  \right]  _{\infty
}\subset\left[  \overline{\operatorname*{co}}\left(  \cup_{t\in T_{1}}%
\partial_{\hat{\delta}}(\varepsilon_{t}f_{t})(x)\right)  \right]  _{\infty},
\]
and we are done thanks to Theorem \ref{p1}.
\end{dem}

The following corollary provides another representation of $\mathrm{N}%
_{\operatorname*{dom}f}(x),$ using the positive part of the $f_{t}$'s instead
of the weighted functions $\varepsilon_{t}f_{t},$ $t\in T\setminus
T_{\varepsilon}(x).$

\begin{cor}
\label{normalnew}Given\textbf{\ }$x\in\operatorname*{dom}f\mathbf{\ }$and
$\varepsilon>0,$ under hypothesis \emph{(\ref{SH})} we have that\textbf{ }%
\begin{equation}
\mathrm{N}_{\operatorname*{dom}f}(x)=\left[  \overline{\operatorname*{co}%
}\left(  \left(
{\textstyle\bigcup\nolimits_{t\in T_{\varepsilon}(x)}}
\partial_{\varepsilon}f_{t}(x)\right)  \cup\left(
{\textstyle\bigcup\nolimits_{t\in T\setminus T_{\varepsilon}(x)}}
\partial_{\varepsilon}f_{t}^{+}(x)\right)  \right)  \right]  _{\infty},
\label{mmain}%
\end{equation}
consequently,
\begin{equation}
\mathrm{N}_{\operatorname*{dom}f}(x)=\left[  \overline{\operatorname*{co}%
}\left(  \left(
{\textstyle\bigcup\limits_{t\in T_{\varepsilon}^{+}(x)}}
\partial_{\varepsilon}f_{t}(x)\right)  \cup\left(
{\textstyle\bigcup\limits_{\substack{t\in T\setminus T_{\varepsilon}%
^{+}(x)\\0\leq\lambda\leq1}}}
\partial_{\varepsilon+\lambda f_{t}(x)}(\lambda f_{t})(x)\right)  \right)
\right]  _{\infty}, \label{kh}%
\end{equation}
where
\[
T_{\varepsilon}^{+}(x):=T_{\varepsilon}(x)\cup\left\{  t\in T:\text{ }%
f_{t}(x)\geq0\right\}  .
\]

\end{cor}

\begin{dem}
We may assume that $f(x)=0.$ First, observe that
\[
\left\{  f_{t},t\in T_{\varepsilon}^{+}(x);\text{ }f_{t}^{+},t\in T\setminus
T_{\varepsilon}^{+}(x)\right\}  \subset\Gamma_{0}(X)
\]
and satisfies
\[
\min\left\{  \inf_{t\in T_{\varepsilon}^{+}(x)}f_{t}(x),\text{ }\inf_{t\in
T\setminus T_{\varepsilon}^{+}(x)}f_{t}^{+}(x)\right\}  \geq-\varepsilon;
\]
that is, (\ref{cb2}) follows and, so, the first statement in Theorem \ref{p1}
implies that
\[
\left[  \overline{\operatorname*{co}}\left(  \left(
{\textstyle\bigcup\nolimits_{t\in T_{\varepsilon}^{+}(x)}}
\partial_{\varepsilon}f_{t}(x)\right)  \cup\left(
{\textstyle\bigcup\nolimits_{t\in T\setminus T_{\varepsilon}^{+}(x)}}
\partial_{\varepsilon}f_{t}^{+}(x)\right)  \right)  \right]  _{\infty}%
\subset\mathrm{N}_{\operatorname*{dom}f}(x).
\]

To establish the opposite inclusion we argue as in the proof of the second
statement in Theorem \ref{p1}. We introduce the nonempty set
\[
D_{\varepsilon}:=\left(
{\textstyle\bigcup\nolimits_{t\in T_{\varepsilon}^{+}(x)}}
\partial_{\varepsilon}f_{t}(x)\right)  \cup\left(
{\textstyle\bigcup\nolimits_{t\in T\setminus T_{\varepsilon}^{+}(x)}}
\partial_{\varepsilon}f_{t}^{+}(x)\right)  ,
\]
and proceed by showing that
\[
\operatorname*{cl}(\operatorname*{dom}\mathrm{\sigma}_{D_{\varepsilon}%
})\subset\operatorname*{cl}(\mathbb{R}_{+}(\operatorname*{dom}f-x)).
\]
Take
\begin{align*}
z\in\operatorname*{dom}\mathrm{\sigma}_{D_{\varepsilon}}  &
=\operatorname*{dom}\left(  \max\left\{  \sup_{t\in T_{\varepsilon}^{+}%
(x)}\mathrm{\sigma}_{\partial_{\varepsilon}f_{t}(x)},\sup_{t\in T\setminus
T_{\varepsilon}^{+}(x)}\mathrm{\sigma}_{\partial_{\varepsilon}f_{t}^{+}%
(x)}\right\}  \right) \\
&  =\operatorname*{dom}\left(  \max\left\{  \sup_{t\in T_{\varepsilon}^{+}%
(x)}(f_{t})_{\varepsilon}^{\prime}(x;\cdot),\sup_{t\in T\setminus
T_{\varepsilon}^{+}(x)}(f_{t}^{+})_{\varepsilon}^{\prime}(x;\cdot)\right\}
\right)  .
\end{align*}
Then
\begin{align*}
z  &  \in\left(  \cap_{t\in T_{\varepsilon}^{+}(x)}\operatorname*{dom}%
((f_{t})_{\varepsilon}^{\prime}(x;\cdot))\right)  \cap\left(  \cap_{T\setminus
T_{\varepsilon}^{+}(x)}\operatorname*{dom}((f_{t}^{+})_{\varepsilon}^{\prime
}(x;\cdot))\right) \\
&  =\left(  \cap_{t\in T_{\varepsilon}^{+}(x)}\mathbb{R}_{+}%
(\operatorname*{dom}f_{t}-x)\right)  \cap\left(  \cap_{T\setminus
T_{\varepsilon}^{+}(x)}\mathbb{R}_{+}(\operatorname*{dom}f_{t}^{+}-x)\right)
\\
&  =\left(  \cap_{t\in T_{\varepsilon}^{+}(x)}\mathbb{R}_{+}%
(\operatorname*{dom}f_{t}-x)\right)  \cap\left(  \cap_{T\setminus
T_{\varepsilon}^{+}(x)}\mathbb{R}_{+}(\operatorname*{dom}f_{t}-x)\right) \\
&  =\cap_{t\in T}\mathbb{R}_{+}(\operatorname*{dom}f_{t}-x),
\end{align*}
and (\ref{egconos}) gives rise to
\[
z\in\cap_{t\in T}\mathbb{R}_{+}(\operatorname*{dom}f_{t}-x)=\mathbb{R}%
_{+}(\operatorname*{dom}f-x)\subset\operatorname*{cl}(\mathbb{R}%
_{+}(\operatorname*{dom}f-x)),
\]
as required.

Finally, the last statement of the corollary follows\ from Lemma \ref{hamu3l}.
\end{dem}

\begin{rem}
Let us note that formula (\ref{kh}) can be simplified, observing for instance
that for all $t\in T\setminus T_{\varepsilon}^{+}(x),$
\[
\partial_{\varepsilon}(0f_{t})(x)=\mathrm{N}_{\operatorname*{dom}f_{t}%
}^{\varepsilon}(x),
\]%
\[
\partial_{\varepsilon+\lambda f_{t}(x)}(\lambda f_{t})(x)=\lambda
\partial_{\lambda^{-1}\varepsilon+f_{t}(x)}f_{t}(x),\text{ for all }%
0<\lambda\leq1.
\]

\end{rem}

\section{Alternative representations of the subdifferential\label{secsub}}

We are dealing again with a nonempty family $\{f_{t},$ $t\in T\}\subset
\Gamma_{0}(X)$ and its\textbf{ }supremum function $f:=\sup_{t\in T}f_{t}.$ As
we remembered in the introduction, in the general case, when no assumption is
made neither on $T$ nor on the mappings $t\mapsto f_{t}(z),$ $z\in X,$ the
subdifferential of $f$ at a point $x\in\operatorname*{dom}f$ is given by
(\cite{HLZ08})\textbf{ }%
\begin{equation}
\partial f(x)=%
{\textstyle\bigcap\nolimits_{L\in\mathcal{F}(x),\text{ }\varepsilon>0}}
\overline{\operatorname*{co}}\left(
{\textstyle\bigcup\nolimits_{t\in T_{\varepsilon}(x)}}
\partial_{\varepsilon}f_{t}(x)+\mathrm{N}_{L\cap\operatorname*{dom}%
f}(x)\right)  \mathbf{,} \label{mfr2}%
\end{equation}
where
\[
\mathcal{F}(x):=\left\{  L\subset X:\text{ }L\text{\ finite-dimensional
subspace with }x\in L\right\}  .
\]

We consider in this section the same family $\left\{  f_{t},t\in T\right\}  $
satisfying\ the standard hypothesis (\ref{SH}); i.e., $T$ is Hausdorff compact
and the mappings $t\mapsto f_{t}(z),$ $z\in X,$ are usc. In such a case,
instead of (\ref{mfr2}) we have the following more precise characterization of
the subdifferential of $f$ (see \cite[Theorem 3.8]{CHL19b}),%
\begin{equation}
\partial f(x)=%
{\textstyle\bigcap\nolimits_{L\in\mathcal{F}(x),\text{ }\varepsilon>0}}
\overline{\operatorname*{co}}\left(
{\textstyle\bigcup\nolimits_{t\in T(x)}}
\partial_{\varepsilon}f_{t}(x)+\mathrm{N}_{L\cap\operatorname*{dom}%
f}(x)\right)  , \label{mfr2b}%
\end{equation}
where we use the active set $T(x)$ instead of $T_{\varepsilon}(x).$

Our objective in this section is to give alternative representations to
(\ref{mfr2b}) for $\partial f(x),$ which are free of $\mathrm{N}%
_{L\cap\operatorname*{dom}f}(x)$, $L\in\mathcal{F}(x).$ The main tools will be
the characterizations of the normal cone to $\operatorname*{dom}f$
provided\ in\ the previous section.

The general characterization of $\partial f(x)$ is given in Theorem
\ref{t1bis}, but we prefer to establish first a preliminary version of it,
which is valid\ in the relevant case when $f$ attains its minimum at\ $x.$

\begin{theo}
\label{t1} Assume that hypothesis \emph{(\ref{SH}) }fulfills. Consider
$x\in\operatorname*{dom}f$ and let $0<\rho_{t}\leq1,$ $t\in T,$ be such that
\[
\inf_{t\in T}(\rho_{t}f_{t})(x)>-\infty.
\]
Then we have that
\begin{equation}
\partial f(x)\subset%
{\textstyle\bigcap\nolimits_{\varepsilon>0}}
\overline{\operatorname*{co}}\left(  \left(
{\textstyle\bigcup\nolimits_{t\in T(x)}}
\partial_{\varepsilon}f_{t}(x)\right)  \cup\left(
{\textstyle\bigcup\nolimits_{t\in T\setminus T(x)}}
\varepsilon\partial_{\varepsilon}(\rho_{t}f_{t})(x)\right)  \right)  .
\label{afms1}%
\end{equation}
Moreover, if $f$ attains its minimum at\ $x,$ then \emph{(\ref{afms1})}
becomes an equality.
\end{theo}

\begin{dem}
Fix $x\in\operatorname*{dom}f$ and assume, without loss of generality, that
$f(x)=0.$ Fix $\varepsilon>0,$ $U\in\mathcal{N},$ and pick\ $L\in
\mathcal{F}(x)$ such that $L^{\bot}\subset U.$ Observe that the family
$\left\{  f_{t},t\in T;\text{ }\mathrm{I}_{L}\right\}  \subset\Gamma_{0}(X)$
also satisfies hypothesis (\ref{SH}) as we can assign to the function
$\mathrm{I}_{L}$ an (isolated) index not belonging to $T.$ Therefore, by
applying Theorem \ref{p1} to the family $\left\{  f_{t},t\in T;\text{
}\mathrm{I}_{L}\right\}  ,$ we obtain that
\begin{align*}
\mathrm{N}_{L\cap\operatorname*{dom}f}(x)  &  =\left[  \overline
{\operatorname*{co}}\left(  \left(
{\textstyle\bigcup\nolimits_{t\in T(x)}}
\partial_{\varepsilon}f_{t}(x)\cup L^{\bot}\right)  \cup\left(  \left(
{\textstyle\bigcup\nolimits_{t\in T\setminus T(x)}}
\partial_{\varepsilon}(\rho_{t}f_{t})(x)\right)  \right)  \right)  \right]
_{\infty}\\
&  =\left[  \overline{\operatorname*{co}}\left(  \left(
{\textstyle\bigcup\nolimits_{t\in T(x)}}
\partial_{\varepsilon}f_{t}(x)\cup L^{\bot}\right)  \cup\left(
{\textstyle\bigcup\nolimits_{t\in T\setminus T(x)}}
\varepsilon\partial_{\varepsilon}(\rho_{t}f_{t})(x)\right)  \right)  \right]
_{\infty},
\end{align*}
where the last equality is a consequence of Lemma \ref{lemconsum0}. Moreover,
applying Lemma \ref{consum-1}, we have that
\begin{equation}
\mathrm{N}_{L\cap\operatorname*{dom}f}(x)=\left[  \overline{\operatorname*{co}%
}\left(  \left(
{\textstyle\bigcup\nolimits_{t\in T(x)}}
\partial_{\varepsilon}f_{t}(x)\right)  \cup\left(
{\textstyle\bigcup\nolimits_{t\in T\setminus T(x)}}
\varepsilon\partial_{\varepsilon}(\rho_{t}f_{t})(x)+L^{\bot}\right)  \right)
\right]  _{\infty}. \label{ar1}%
\end{equation}
Next, by combining this relation\ and\ (\ref{mfr2b}),
\begin{align*}
\partial f(x)  &  \subset\overline{\operatorname*{co}}\left(
{\textstyle\bigcup\nolimits_{t\in T(x)}}
\partial_{\varepsilon}f_{t}(x)+\mathrm{N}_{L\cap\operatorname*{dom}%
f}(x)\right) \\
&  =\overline{\operatorname*{co}}\left(
{\textstyle\bigcup\limits_{t\in T(x)}}
\partial_{\varepsilon}f_{t}(x)+\left[  \overline{\operatorname*{co}}\left(
\left(
{\textstyle\bigcup\limits_{t\in T(x)}}
\partial_{\varepsilon}f_{t}(x)\right)  \cup\left(
{\textstyle\bigcup\limits_{t\in T\setminus T(x)}}
\varepsilon\partial_{\varepsilon}(\rho_{t}f_{t})(x)+L^{\bot}\right)  \right)
\right]  _{\infty}\right) \\
&  \subset\overline{\operatorname*{co}}\left(  \left(
{\textstyle\bigcup\nolimits_{t\in T(x)}}
\partial_{\varepsilon}f_{t}(x)\right)  \cup\left(
{\textstyle\bigcup\nolimits_{t\in T\setminus T(x)}}
\varepsilon\partial_{\varepsilon}(\rho_{t}f_{t})(x)+L^{\bot}\right)  \right)
\\
&  \subset\operatorname*{co}\left(  \left(
{\textstyle\bigcup\nolimits_{t\in T(x)}}
\partial_{\varepsilon}f_{t}(x)\right)  \cup\left(
{\textstyle\bigcup\nolimits_{t\in T\setminus T(x)}}
\varepsilon\partial_{\varepsilon}(\varepsilon_{t}f_{t})(x)\right)  \right)
+2U.
\end{align*}
Consequently, (\ref{afms1}) follows by intersecting over $\varepsilon>0$ and
$U\in\mathcal{N}.$

We proceed now by showing the opposite inclusion in (\ref{afms1}) when $x$ is
a minimizer of $f.$ By the current assumption, we choose an $M>0$ such that
\[
\inf_{t\in T}(\rho_{t}f_{t})(x)\geq-M,
\]
and take $x^{\ast}$ in the right-hand side of (\ref{afms1}); that is, for each
fixed $\varepsilon>0,$
\begin{equation}
x^{\ast}\in\overline{\operatorname*{co}}\left(  \left(
{\textstyle\bigcup\nolimits_{t\in T(x)}}
\partial_{\varepsilon}f_{t}(x)\right)  \cup\left(
{\textstyle\bigcup\nolimits_{t\in T\setminus T(x)}}
\varepsilon\partial_{\varepsilon}(\rho_{t}f_{t})(x)\right)  \right)  .
\label{doi}%
\end{equation}
Observe that, if $z^{\ast}\in\partial_{\varepsilon}f_{t}(x),$ $t\in T(x),$
then for every $z\in X$
\[
\left\langle z^{\ast},z-x\right\rangle \leq f_{t}(z)-f_{t}(x)+\varepsilon\leq
f(z)+\varepsilon,
\]
and so $z^{\ast}\in\partial_{\varepsilon}f(x).$ Moreover, if $z^{\ast}%
\in\partial_{\varepsilon}(\rho_{t}f_{t})(x),$ $t\in T\setminus T(x),$ then for
every $z\in X$
\begin{align*}
\left\langle z^{\ast},z-x\right\rangle  &  \leq\rho_{t}f_{t}(z)-\rho_{t}%
f_{t}(x)+\varepsilon\\
&  \leq\rho_{t}f(z)-M+\varepsilon\leq f(z)+M+\varepsilon,
\end{align*}
since that $x$ is a minimizer of $f$ and $f(z)\geq f(x)=0.$ Hence, $z^{\ast
}\in\partial_{\varepsilon+M}f_{t}(x).$ Consequently, taking into account Lemma
\ref{lemvo} and that $\varepsilon>0$ was arbitrarily chosen, (\ref{doi}) leads
us\ to
\[
x^{\ast}\in\cap_{\varepsilon>0}\overline{\operatorname*{co}}\left(
\partial_{\varepsilon}f(x)\cup\varepsilon\partial_{\varepsilon+M}f(x)\right)
=\partial f(x),
\]
that is, the opposite inclusion in (\ref{afms1}) is also true.
\end{dem}

\begin{exam}
\emph{The inclusion in (\ref{afms1}) may be strict when }$x$\emph{ is not a
minimizer of }$f.$\emph{ Consider the family }%
\[
\left\{  f_{t},\text{ }t\in T;\text{ }h\right\}  ,
\]
\emph{where }$h$\emph{ is the constant function }$h\equiv f(x)-1,$\emph{ and
denote by }$g$\emph{ the associated supremum function.} \emph{Then, by the
lower semicontinuity of the }$f_{t}$\emph{'s, the functions }$f$\emph{ and
}$g$\emph{ coincide in a neighborhood of }$x$\emph{, entailing\ }%
\[
\partial f(x)=\partial g(x).
\]
\emph{If (\ref{afms1}) would be an equality, with any weighting parameter
associated to\ }$h,$ \emph{then by taking into account that }$\partial
_{\varepsilon}h(x)=\left\{  \theta\right\}  $\emph{ we would have }%
\[
\partial g(x)=%
{\textstyle\bigcap\nolimits_{\varepsilon>0}}
\overline{\operatorname*{co}}\left(  \left(
{\textstyle\bigcup\nolimits_{t\in T(x)}}
\partial_{\varepsilon}f_{t}(x)\right)  \cup\left(
{\textstyle\bigcup\nolimits_{t\in T\setminus T(x)}}
\varepsilon\partial_{\varepsilon}(\rho_{t}f_{t})(x)%
{\textstyle\bigcup}
\left\{  \theta\right\}  \right)  \right)  ,
\]
\emph{but this implies\ that }%
\[
\theta\in\partial g(x)=\partial f(x),
\]
\emph{which contradicts our assumption that }$x$\emph{ is not a minimizer of
}$f.$
\end{exam}

\begin{theo}
\label{t1bis} Assume that hypothesis \emph{(\ref{SH}) }fulfills. Consider
$x\in\operatorname*{dom}f$ and let $0<\rho_{t}\leq1,$ $t\in T,$ be such that
\[
\inf_{t\in T}(\rho_{t}f_{t})(x)>-\infty.
\]
Then we have that
\begin{equation}
\partial f(x)=%
{\textstyle\bigcap\nolimits_{\varepsilon>0}}
\overline{\operatorname*{co}}\left(  \left(
{\textstyle\bigcup\nolimits_{t\in T(x)}}
\partial_{\varepsilon}f_{t}(x)\right)  +\left(
{\textstyle\bigcup\nolimits_{t\in T\setminus T(x)}}
\left\{  0,\varepsilon\right\}  \partial_{\varepsilon}(\rho_{t}f_{t}%
)(x)\right)  \right)  \label{afms2}%
\end{equation}
(with $%
{\textstyle\bigcup\nolimits_{\emptyset}}
=\left\{  \theta\right\}  $ when $T(x)=T$).
\end{theo}

\begin{rem}
[before the proof]\emph{Observe that the operator} $\operatorname*{co}$
\emph{determines that (\ref{afms2}) can be equivalently written as}
\[
\partial f(x)=%
{\textstyle\bigcap\nolimits_{\varepsilon>0}}
\overline{\operatorname*{co}}\left(  \left(
{\textstyle\bigcup\nolimits_{t\in T(x)}}
\partial_{\varepsilon}f_{t}(x)\right)  +\left(
{\textstyle\bigcup\nolimits_{t\in T\setminus T(x)}}
\left[  0,\varepsilon\right]  \partial_{\varepsilon}(\rho_{t}f_{t})(x)\right)
\right)  .
\]

\end{rem}

\begin{dem}
Take $x\in\operatorname*{dom}f$ such that $f(x)=0$ (without loss of
generality). Fix $\varepsilon>0,$ $U\in\mathcal{N},$ and pick\ $L\in
\mathcal{F}(x)$ such that $L^{\bot}\subset U.$ By arguing\ as in the beginning
of the proof of (\ref{afms1}), and taking into account\ Lemma \ref{consum-1},
we obtain that
\begin{align}
\mathrm{N}_{L\cap\operatorname*{dom}f}(x) &  =\left[  \overline
{\operatorname*{co}}\left(  \left(
{\textstyle\bigcup\nolimits_{t\in T(x)}}
\partial_{\varepsilon}f_{t}(x)\right)  +\left(
{\textstyle\bigcup\nolimits_{t\in T\setminus T(x)}}
\varepsilon\partial_{\varepsilon}(\rho_{t}f_{t})(x)\right)  +L^{\bot}\right)
\right]  _{\infty}\label{ar}\\
&  \subset\left[  \overline{\operatorname*{co}}\left(  \left(
{\textstyle\bigcup\nolimits_{t\in T(x)}}
\partial_{\varepsilon}f_{t}(x)\right)  +\left(
{\textstyle\bigcup\nolimits_{t\in T\setminus T(x)}}
\left\{  0,\varepsilon\right\}  \partial_{\varepsilon}(\rho_{t}f_{t}%
)(x)\right)  +L^{\bot}\right)  \right]  _{\infty}.\nonumber
\end{align}
(Observe the difference between (\ref{ar}) and (\ref{ar1}).)

Due to the lower semicontinuity of the $(\rho_{t}f_{t})$'s, the sets
$\partial_{\varepsilon}(\rho_{t}f_{t})(x)$ are nonempty and we have
\[
\theta\in%
{\textstyle\bigcup\nolimits_{t\in T\setminus T(x)}}
\left\{  0,\varepsilon\right\}  \partial_{\varepsilon}(\rho_{t}f_{t}%
)(x)+L^{\bot}.
\]
Thus,\ using again (\ref{mfr2b}),%
\begin{align*}
\partial f(x) &  \subset\overline{\operatorname*{co}}\left(
{\textstyle\bigcup\nolimits_{t\in T(x)}}
\partial_{\varepsilon}f_{t}(x)+\mathrm{N}_{L\cap\operatorname*{dom}%
f}(x)\right)  \\
&  \subset\overline{\operatorname*{co}}\left(
{\textstyle\bigcup\nolimits_{t\in T(x)}}
\partial_{\varepsilon}f_{t}(x)+%
{\textstyle\bigcup\nolimits_{t\in T\setminus T(x)}}
\left\{  0,\varepsilon\right\}  \partial_{\varepsilon}(\rho_{t}f_{t}%
)(x)+L^{\bot}\right)  \\
&  \subset\operatorname*{co}\left(
{\textstyle\bigcup\nolimits_{t\in T(x)}}
\partial_{\varepsilon}f_{t}(x)+%
{\textstyle\bigcup\nolimits_{t\in T\setminus T(x)}}
\left\{  0,\varepsilon\right\}  \partial_{\varepsilon}(\rho_{t}f_{t}%
)(x)\right)  +2U.
\end{align*}
Therefore the first inclusion \textquotedblleft$\subset$\textquotedblright\ in
(\ref{afms2}) follows by intersecting over $\varepsilon>0$ and $U\in
\mathcal{N}.$

Conversely, to show the inclusion \textquotedblleft$\supset$\textquotedblright%
\ in (\ref{afms2}), we take $x^{\ast}$ in the right-hand side of (\ref{afms2})
and choose an $M>0$ such that
\[
\inf_{t\in T}(\rho_{t}f)(x)>-M.
\]
Thus, for each $\varepsilon>0$ we obtain (similarly\ to the last part of the
proof of Theorem \ref{t1})%
\begin{align}
x^{\ast} &  \in\overline{\operatorname*{co}}\left(  \left(
{\textstyle\bigcup\nolimits_{t\in T(x)}}
\partial_{\varepsilon}f_{t}(x)\right)  +\left(
{\textstyle\bigcup\nolimits_{t\in T\setminus T(x)}}
\left\{  0,\varepsilon\right\}  \partial_{\varepsilon}(\rho_{t}f_{t}%
)(x)\right)  \right)  \nonumber\\
&  \subset\overline{\operatorname*{co}}\left(  \partial_{\varepsilon
}f(x)+\left(
{\textstyle\bigcup\nolimits_{t\in T\setminus T(x)}}
\left\{  0,\varepsilon\right\}  \partial_{M+\varepsilon}(\rho_{t}f)(x)\right)
\right)  \nonumber
\end{align}
(Observe that, compared to the proof of Theorem \ref{t1}, here we maintain the
parameters $\rho_{t},$ $t\in T\setminus T(x),$ because $x$ is not necessarily
a minimizer of $f.$)

Writing, due to last relation,
\begin{align*}
x^{\ast} &  \in\operatorname*{cl}\left(  \operatorname*{co}\left(
\partial_{\varepsilon}f(x)+%
{\textstyle\bigcup\nolimits_{t\in T\setminus T(x)}}
\left\{  0,\varepsilon\right\}  \partial_{M+\varepsilon}(\rho_{t}f)(x)\right)
\right)  \\
&  =\operatorname*{cl}\left(  \partial_{\varepsilon}f(x)+\operatorname*{co}%
\left(
{\textstyle\bigcup\nolimits_{t\in T\setminus T(x)}}
\left\{  0,\varepsilon\right\}  \partial_{M+\varepsilon}(\rho_{t}f)(x)\right)
\right)  ,
\end{align*}
leading to the existence of\ nets $(y_{i}^{\ast})_{i}\subset\partial
_{\varepsilon}f(x),$ $(\lambda_{i,k})_{i}\subset\left[  0,1\right]  ,$
$\Sigma_{k=1,\cdots,k_{i}}\lambda_{i,k}\leq1,$ $(t_{i,k})_{i}\subset
T\setminus T(x),$ and $(z_{i,k}^{\ast})_{i}\subset\partial_{M+\varepsilon
}(\rho_{t_{i,k}}f)(x),$ $k=1,\cdots,k_{i},$ $k_{i}\geq1,$ such that
\[
x^{\ast}=\lim_{i}\left(  y_{i}^{\ast}+%
{\textstyle\sum\limits_{k=1,\cdots,k_{i}}}
\varepsilon\lambda_{i,k}z_{i,k}^{\ast}\right)  .
\]
Consequently, since $f(x)=0,$ for each\ $y\in\operatorname*{dom}f$%
\begin{align*}
\left\langle x^{\ast},y-x\right\rangle  &  =\lim_{i}\left\langle y_{i}^{\ast
}+\varepsilon%
{\textstyle\sum\limits_{k=1,\cdots,k_{i}}}
\lambda_{i,k}z_{i,k}^{\ast},y-x\right\rangle \\
&  \leq\limsup_{i}\left(  \left(  f(y)-f(x)+\varepsilon\right)  +\varepsilon
\left(
{\textstyle\sum\limits_{k=1,\cdots,k_{i}}}
\lambda_{i,k}(\rho_{t_{i,k}}f(y)-\rho_{t_{i,k}}f(x)+M+\varepsilon)\right)
\right)  \\
&  =\limsup_{i}\left(  f(y)+\varepsilon+\varepsilon%
{\textstyle\sum\limits_{k=1,\cdots,k_{i}}}
\lambda_{i,k}(\rho_{t_{i,k}}f(y)+M+\varepsilon)\right)  \\
&  \leq f(y)+\varepsilon+\varepsilon(f^{+}(y)+M+\varepsilon),
\end{align*}
and\textbf{ }$x^{\ast}\in\partial_{\varepsilon}f(x),$ by taking\ $\varepsilon
\downarrow0.$
\end{dem}

\begin{cor}
\label{corr} Assume that hypothesis \emph{(\ref{SH}) }fulfills. Then for every
$x\in\operatorname*{dom}f$ we have that
\[
\partial f(x)=%
{\textstyle\bigcap\nolimits_{\varepsilon>0}}
\overline{\operatorname*{co}}\left(  \left(
{\textstyle\bigcup\nolimits_{t\in T(x)}}
\partial_{\varepsilon}f_{t}(x)\right)  +\left(
{\textstyle\bigcup\nolimits_{t\in T\setminus T(x)}}
\left\{  0,\varepsilon\right\}  \partial_{\varepsilon}(\rho_{t}f_{t}%
)(x)\right)  \right)  ,
\]
where $\rho_{t}=\rho_{t}(\varepsilon)$ is defined as
\[
\rho_{t}:=\frac{\varepsilon}{2f(x)-2f_{t}(x)+\varepsilon},\text{ \ \ }t\in
T\setminus T(x).
\]
In particular, if $f$ attains its minimum at $x,$ then we also have that
\[
\partial f(x)=%
{\textstyle\bigcap\nolimits_{\varepsilon>0}}
\overline{\operatorname*{co}}\left(  \left(
{\textstyle\bigcup\nolimits_{t\in T(x)}}
\partial_{\varepsilon}f_{t}(x)\right)  \cup\left(
{\textstyle\bigcup\nolimits_{t\in T\setminus T(x)}}
\varepsilon\partial_{\varepsilon}(\rho_{t}f_{t})(x)\right)  \right)  .
\]

\end{cor}

\begin{dem}
It suffices to apply Theorems \ref{t1} and \ref{t1bis} by replacing the
parameters $\rho_{t},$ $t\in T,$ there by
\[
\hat{\rho}_{t}:=1,\text{ if }t\in T(x),\text{ }\hat{\rho}_{t}:=\frac
{\varepsilon}{2f(x)-2f_{t}(x)+\varepsilon}\text{, }t\in T\setminus T(x).
\]
Indeed, for all $t\in T$ we have that $0<\rho_{t}<1$ and
\[
\hat{\rho}_{t}f_{t}(x)\geq\min\left\{  0,\inf_{t\in T\setminus T(x)}%
\frac{\varepsilon f_{t}(x)}{2f(x)-2f_{t}(x)+\varepsilon}\right\}  \geq
-\frac{\varepsilon}{2}.
\]

\end{dem}

\begin{cor}
Assume that \emph{(\ref{SH})} fulfills.\ If $x\in\operatorname*{dom}f$ is such
that
\[
\inf_{t\in T}f_{t}(x)>-\infty,
\]
then we have
\[
\partial f(x)=%
{\textstyle\bigcap\nolimits_{\varepsilon>0}}
\overline{\operatorname*{co}}\left(  \left(
{\textstyle\bigcup\nolimits_{t\in T(x)}}
\partial_{\varepsilon}f_{t}(x)\right)  +\left(
{\textstyle\bigcup\nolimits_{t\in T\setminus T(x)}}
\left\{  0,\varepsilon\right\}  \partial_{\varepsilon}f_{t}(x)\right)
\right)  ,
\]
and, when\ additionally $f$ attains its minimum at $x,$%
\[
\partial f(x)=%
{\textstyle\bigcap\nolimits_{\varepsilon>0}}
\overline{\operatorname*{co}}\left(  \left(
{\textstyle\bigcup\nolimits_{t\in T(x)}}
\partial_{\varepsilon}f_{t}(x)\right)  \cup\left(
{\textstyle\bigcup\nolimits_{t\in T\setminus T(x)}}
\varepsilon\partial_{\varepsilon}f_{t}(x)\right)  \right)  .
\]

\end{cor}

An obvious\ consequence of Theorem \ref{t1bis} is the following extension
of\ the Brøndsted formula in \cite{Br72}, and the formula given in
\cite[Proposition 6.3]{HL08} (in finite dimensions and under the continuity of
the $f_{t}$'s).

\begin{cor}
Assume that \emph{(\ref{SH})} fulfills.\ If $x\in\operatorname*{dom}f$ is such
that $T(x)=T,$ then\
\[
\partial f(x)=%
{\textstyle\bigcap\nolimits_{\varepsilon>0}}
\overline{\operatorname*{co}}\left(
{\textstyle\bigcup\nolimits_{t\in T(x)}}
\partial_{\varepsilon}f_{t}(x)\right)  .
\]

\end{cor}

\begin{dem}
It is immediate from Theorem \ref{t1bis}, since that $\inf_{t\in T}%
f_{t}(x)=\inf_{t\in T(x)}f_{t}(x)=0.$
\end{dem}

We close the paper by deriving new optimality conditions for the following
convex optimization problem with infinitely many constraints,
\[
\mathcal{(P)}:\text{ \ }\operatorname*{Inf}g(x),\ \text{ subject to }%
f_{t}(x)\leq0,\text{ }t\in T,
\]
where $T$ is an arbitrary (possibly, infinite) set. We refer e.g. to
\cite{GL98}, \cite{GL18}, \cite{LoSt07}, etc., and references therein, for
theory, algorithms and applications of this model. We have the following
result in the continuous framework; i.e., $T$ is a Hausdorff compact
set$\ $and the family\ $\left\{  g;\text{ }f_{t},\text{ }t\in T\right\}
\subset\Gamma_{0}(X)$ satisfies\ condition (\ref{SH})\textbf{.} See,
also,\ \cite{CHL16, CHL20} for optimality conditions for $\mathcal{(P)}$ in
different frameworks.

\begin{cor}
Let $\bar{x}\ $be an optimal\ solution\ of $\mathcal{(P)}$, and take
$0<\rho_{t}\leq1,$ $t\in T\setminus T(\bar{x}),$ such that
\[
\inf_{t\in T\setminus T(\bar{x})}(\rho_{t}f_{t})(\bar{x})>-\infty.
\]
Then, for every $\varepsilon>0$ and every $U\in\mathcal{N}$, there are
associated $t_{1},\cdots,t_{m}\in T(\bar{x}),$ $t_{m+1},\cdots,t_{m+n}\in
T\setminus T(\bar{x}),$ and $(\lambda_{0},\lambda_{1},\cdots,\lambda
_{m},\lambda_{m+1},\cdots,\lambda_{m+n})\in\Delta_{m+n+1},$ $m,n\geq1,$ such
that
\[
0_{n}\in\lambda_{0}\partial_{\varepsilon}g(\bar{x})+%
{\textstyle\sum\limits_{i=1}^{m}}
\lambda_{i}\partial_{\varepsilon}f_{t_{i}}(\bar{x})+%
{\textstyle\sum\limits_{i=m+1}^{m+n}}
\varepsilon\lambda_{i}\partial_{\varepsilon}(\rho_{t_{i}}f_{t_{i}})+U.
\]
Moreover\textbf{, }$\lambda_{0}>0$\textbf{ }when the Slater condition is
satisfied; that is, there is some $x_{0}\in X$ such that
\[
f_{t}(x_{0})<0\mathbf{\ \ }\text{for all }t\in T.
\]

\end{cor}

\begin{rem}
[before the proof]\emph{Observe that some of the multipliers} $\lambda
_{0},\lambda_{1},\cdots,\lambda_{m},$ $\lambda_{m+1},\cdots,\lambda_{m+n}$
\emph{can be zero, but their sum is one. Note that, due to the hypothesis
(\ref{SH}), }$T(\bar{x})\neq\emptyset$ \emph{but} $T\setminus T(\bar{x})$
\emph{can be empty. In the last case, the last relation collapses to}
\[
0_{n}\in\lambda_{0}\partial_{\varepsilon}g(\bar{x})+%
{\textstyle\sum\limits_{i=1}^{m}}
\lambda_{i}\partial_{\varepsilon}f_{t_{i}}(\bar{x})+U.
\]

\end{rem}

\begin{dem}
It is easy to see that $\bar{x}$ is a global minimum of the supremum function
$f:X\rightarrow\mathbb{R\cup\{+\infty\}}$, defined as
\[
f(x):=\sup\{g(x)-g(\bar{x}),\ f_{t}(x),t\in T\};
\]
that is, $0_{n}\in\partial f(\bar{x}).$ Then, by Theorem \ref{t1},
\[
0_{n}\in\partial f(\bar{x})=%
{\textstyle\bigcap\nolimits_{\varepsilon>0}}
\overline{\operatorname*{co}}\left(  \left(  \partial_{\varepsilon}g(\bar
{x})\cup\left(
{\textstyle\bigcup\nolimits_{t\in T(\bar{x})}}
\partial_{\varepsilon}f_{t}(\bar{x})\right)  \right)  \cup\left(
{\textstyle\bigcup\nolimits_{t\in T\setminus T(\bar{x})}}
\varepsilon\partial_{\varepsilon}(\rho_{t}f_{t})(\bar{x})\right)  \right)  ,
\]
leading us to the conclusion of the first statement of the corollary.

Now, we suppose that the Slater condition holds; that is, due to (\ref{SH}),
the supremum function $h:=\sup_{t\in T}f_{t}$ satisfies
\[
h(x_{0})<0,
\]
for some $x_{0}\in X.$ Let us\ proceed by contradiction, assuming\ that
$\lambda_{0}=0;$ that is, since $g\in\Gamma_{0}(X)$ and so $\partial
_{\varepsilon}g(\bar{x})\neq\emptyset$,%
\begin{equation}
0_{n}\in%
{\textstyle\sum\limits_{i=1}^{m}}
\lambda_{i}\partial_{\varepsilon}f_{t_{i}}(\bar{x})+%
{\textstyle\sum\limits_{i=m+1}^{m+n}}
\varepsilon\lambda_{i}\partial_{\varepsilon}(\rho_{t_{i}}f_{t_{i}})+U.
\label{rea}%
\end{equation}
Observe that
\[
\partial_{\varepsilon}f_{t_{i}}(\bar{x})\subset\partial_{\varepsilon}h(\bar
{x}),\text{ \ }i=1,\cdots,m,
\]
since that $f_{t_{i}}\leq h$ and $f_{t_{i}}(\bar{x})=h(\bar{x})=0,$
$i=1,\cdots,m.$ Moreover, if\ $M>0$ is such that\ $\inf_{t\in T\setminus
T(\bar{x})}(\rho_{t}f_{t})(\bar{x})>-M,$ then for every\ $x^{\ast}\in
\partial_{\varepsilon}(\rho_{t_{i}}f_{t_{i}})(\bar{x}),$ $t_{i}\in T\setminus
T(\bar{x}),$ we have for all $y\in\operatorname*{dom}h$%
\begin{align*}
\left\langle x^{\ast},y-\bar{x}\right\rangle  &  \leq(\rho_{t_{i}}f_{t_{i}%
})(y)-(\rho_{t_{i}}f_{t_{i}})(\bar{x})+\varepsilon\\
&  \leq(\rho_{t_{i}}h)(y)-\inf_{t\in T\setminus T(\bar{x})}(\rho_{t}%
f_{t})(\bar{x})+\varepsilon\\
&  \leq h^{+}(y)+M+\varepsilon;
\end{align*}
that is, $x^{\ast}\in\partial_{M+\varepsilon}h^{+}(\bar{x}).$ Consequently,
(\ref{rea}) reads
\[
0_{n}\in%
{\textstyle\sum\limits_{i=1}^{m}}
\lambda_{i}\partial_{\varepsilon}h(\bar{x})+\varepsilon%
{\textstyle\sum\limits_{i=m+1}^{m+n}}
\lambda_{i}\partial_{M+\varepsilon}h^{+}(\bar{x})+U\subset\operatorname*{co}%
(\partial_{\varepsilon}h(\bar{x}),\varepsilon\partial_{M+\varepsilon}%
h^{+}(\bar{x}))+U,
\]
and so, according to Lemma \ref{lemvo},
\[
0_{n}\in\cap_{\varepsilon>0}\overline{\operatorname*{co}}(\partial
_{\varepsilon}h(\bar{x}),\varepsilon\partial_{M+\varepsilon}h^{+}(\bar
{x}))=\partial h(\bar{x}).
\]
This is a contradiction because $0=h(\bar{x})\leq h(x_{0})<0.$
\end{dem}

\end{document}